\documentclass{article}
\begin{document}
\begin{center}\Large{An Enumerative Function}\end{center}
 \begin{center}MILAN JANJI\'C\end{center}
 \begin{center}Faculty of Natural Sciences and Mathematics, \end{center}
\begin{center} Banja Luka, Republic of Srpska,
Bosnia and Herzegovina.\end{center}
 \begin{center}e-mail: agnus@blic.net\end{center}

\begin{abstract}
We define an enumerative function $F(n,k,P,m)$ which is  a generalization of binomial coefficients.
Special cases of this function are also power function, factorials, rising factorials and  falling factorials.

The first section of the paper is an introduction.

In the second  section we derive an explicit formula for $F.$
From the expression for the power function we obtain a number theory result.

Then we derive a formula which shows that the case of arbitrary $m$ may be reduced to the case $m=0.$
This formula extends  Vandermonde convolution.

In the second section we describe $F$ by the series of recurrence relations
with respect to each of arguments $k,\;n,\;$ and $P.$ As a special case of the first recurrence relation we state a binomial identity.
As a consequence of the second recurrence relation  we obtain  relation for coefficients of Chebyshev polynomial of both kind.
This means that these polynomials might be defined in pure combinatorial way.
\end{abstract}

\vspace{5mm} \noindent Keywords: Enumerative functions, Vandermonde convolution, Chebyshev polynomials of the
first kind, Chebyshev polynomials of the second kind


\vspace{5mm}
\section{Introduction}
For a set  $P=\{p_1,p_2,\ldots,p_n\}$  of positive integers and
nonnegative integers  $m$ and $k,$ we consider  the set $X$ consisting of
$n$ blocks $X_i,\;(|X_i|=p_i,\;i=1,2,\ldots,n),$ which we shall called the main blocks of $X,$ and an additional
$m$-block $X_{n+1}.$ We shall denote
$|X|=N=p_1+\cdots+p_n+m.$

\noindent\textbf{Definition 1.} \textit{
By an $s$-inset of $X$ we shall
mean some $s$-subset of $X$ intersecting each main block of $X.$}

\noindent\textbf{Definition 2.} \textit{ We
define the function $F(n,k,P,m)$ to be the number of $n+k$-insets
of $X.$ We also define $F(n,k,P,m)=0$ if $k<0.$}

The function $F(n,k,P,m)$ seems to be  a natural
 generalization of binomial coefficients, since in the case that  all main blocks have only one elements $F(n,k,P,m)$ becomes the
  binomial coefficient ${m\choose k}$.

Some other combinatorial functions also arise from $F(n,k,P,n).$

In the case that all main blocks have the same number ($p$) of elements then
$$F(n,0,P,m)=p^n.$$

For $p_i=s+i,\;(i=1,2,\ldots,n)$ rising factorials arises.
$$F(n,0,P,m)=(s+1)(s+2)\cdots(s+n),$$

For $p_i=n+s-i,\;(i=1,2,\ldots,n)$ falling factorial are obtained.
$$F(n,0,P,m)=(n+s-1)(n+s-2)\cdots s.$$

In the first section of this paper we prove that $F(n,k,P,m)$ may be expressed as an alternating sum of pairs of binomial coefficients.
This formula takes a simpler form when all main blocks have the same number of elements. It follows that each above mentioned function may be obtained as
such one  alternating sum.
From this fact we obtain a number theory result.

In this section we also prove that $F(n,k,P,m)$ may be expressed
 in terms of $F(n,k_1,P,0),$ that is, in terms of sets with  no additional blocks. The case $n=1$ in thus obtained formula is well-known Vandermonde convolution.

 The second section is devoted to the recurrence relations  of $F(n,k,P,m).$ We prove three formulae, with respect to each parameter $n,\;k,\;P.$
As a consequence of the recurrence relation with respect to $k$ some binomial identities are obtained.

 The recurrence relation with respect to  $n,$ in the case that each main block have exactly  two elements, gives the coefficients of Chebyshev polynomials of both kind. This proves
 that this polynomials may be defined in pure combinatorial way.

We finish the paper with formula in which $F(n,k,P,M)$ is expressed in terms of $F(n_1,k_1,Q,m_1),$ where $q_i=p_i$ or $q_i=p_i-1.$

Note that our function generates a number of sequences in the well-known Sloane's Encyclopedia of Integers Sequences [1].
\section{The formula for $F(n,k,P,m)$}
Note first that it is clear that $U$ is an inset of $X$ if and only if $X\setminus
U$ contain none of main blocks.  Thus, Definition 1
may be stated in the following equivalent form.

\noindent\textbf{Definition 3.} The function $F(n,k,P,m)$ is the
number of $n+k$-subsets $U$ of $X$ such that $X\setminus U$
contain none of the main blocks  $X_i,\;(i=1,2,\ldots,n).$

We shall first derive an explicit formula for $F(n,k,P,m).$

\noindent\textbf{Theorem 1.} \textit{ It holds}
\begin{equation}\label{tau}F(n,k,P,m)=\sum_{I\subseteq [n]}(-1)^{|I|}{N-\sum_{i\in I}p_i\choose
n+k},\end{equation}
where the sum is taken over all subsets of $[n].$

\noindent\textbf{Proof.} For $i=1,2,\ldots,n$
 and a subset $Z$ of $X$ define the property $q_i,\;(i=1,2,\ldots,n)$ to be:
Block $X_i$ does not intersect $Z.$ By PIE method we obtain
$$F(n,k,P,m)=\sum_{I\subseteq [n]}(-1)^{|I|}N(I),$$ where
$N(I)$ is the number of $n+k$-subsets of $X$ which do not
intersect main blocks $X_i,\;(i\in I).$ There are $${N-\sum_{i\in
I}p_i\choose n+k}$$ such subsets and the formula (\ref{tau}) is
proved.

If $p_1=p_2=\cdots=p_n=p$ then the formula (\ref{tau}) takes a
little simpler form
\begin{equation}\label{f1}
F(n,k,P,m)=\sum_{i=0}^n(-1)^i{n\choose i}{np+m-ip\choose
n+k}.\end{equation}

If each of the main blocks has one element, that is, if $p_1=\cdots=p_n=1$ then
$$F(n,k,P,m)={m\choose k},$$
which shows that $F(n,k,P,m)$ is rather natural extension of binomial coefficients.

If $p_1=\cdots=p_n=p>1$ then $$F(n,0,P,m)=p^n,$$ and
\begin{equation}\label{n=1}F(1,k,P,m)={p+m\choose k+1}-{m\choose k+1}.\end{equation}

From this we obtain the following identities
$${m\choose k}=\sum_{i=0}^n(-1)^i{n\choose i}{n+m-i\choose
n+k},$$ which is well-known binomial identity, and

$$p^n=\sum_{i=0}^n(-1)^i{n\choose i}{pn+m-pi\choose
n}.$$

If $n,$ in the preceding equation, is prime, then it divides all terms in the sum on the right side except the first and the last.
The first term is ${n+m\choose n},$ and the last is ${m\choose n}.$ We thus obtain

 \noindent\textbf{Corollary 1.} \textit{ If $q$ is a prime then for integers $r>1$ and $m\geq 0$ holds}
$$q\left|r^q-{rq+m\choose q}+{m\choose q}.\right.$$
Taking $k=0$ in  (\ref{tau}) we obtain
$$F(n,0,P,m)=\prod_{i=1}^np_i.$$
Specially, for $p_i=i,\;(i=1,2,\ldots,n)$ it follows
$$F(n,0,P,m)=n!.$$

The next formula reduces the case of arbitrary $m$ to the case $m=0.$

\noindent\textbf{Theorem 2.} \textit{The following formula holds}
\begin{equation}\label{r2}F(n,k,P,m)=\sum_{i=0}^{min(m,k)}{m\choose i}F(n,k-i,P,0).\end{equation}

\noindent\textbf{Proof.}
 Omitting the $m$-blocks of $X$ we
obtain the set $Y$ with no additional block.

Each $n+k$-inset of $Y$ is an $n+k$-inset of $X.$  There are
$F(n,k,P,0)$ such insets. In this way we obtain all $n+k$-inset of $X$ not intersecting $m$-block.

The remaining $n+k$-insets of $X$ are obtained as a union of  some $n+k-i,\;(1\leq i\leq m)$-inset of $Y$ and some $i$-set of additional block.
In such a way we obtain (\ref{r2}).

If all main blocks has the same number $p$  of elements and if $n=1$ the formula (\ref{r2}) becomes
$$F(1,k,P,m)=\sum_{i=0)}^{min(m,k)}{m\choose i}F(1,k-i,P,0).$$
According to (\ref{n=1}) holds
$$F(1,k-i,P,0)={p\choose k-i+1}.$$

We thus obtain the following identity

$${p+m\choose 1+k}=\sum_{i=0}^{min(m,k+1)}{m\choose i}{p\choose k-i+1},$$
which is well-known Vandermonde convolution.

\section{Recurrence Relations }

The first recurrence relation  is with respect to the parameter $k.$

 \noindent\textbf{Theorem 3.} \textit{For each $s=1,2,\ldots$ the following formula is true}
\begin{equation}\label{r1}F(n,k,P,m)=\sum_{i=0}^s(-1)^i{s\choose i}F(n,k+s,P,m+s-i).\end{equation}

\noindent\textbf{Proof.}
 Taking sum over all $I,\;(\subseteq [n])$ on the left and the right side of  well-known recurrence
 relation for binomial coefficients  $${N-\sum_{i\in I}p_i\choose n+k}+{N-\sum_{i\in I}p_i\choose
n+k+1}={N+1-\sum_{i\in I}p_i\choose n+k+1}$$
we obtain
$$F(n,k,P,m)=F(n,k+1,P,m+1)-F(n,k+1,P,m),$$ and then (\ref{r1}) follows by induction.

In the case $m=0,\;n=1$ the formula (\ref{r2}), for each $s=0,1,2,\ldots,$ produces
the following binomial identities
$${p\choose k+1}=\sum_{i=0}^{s}(-1)^s{s\choose i}{p+s-i\choose k+s+1}.$$

The following recurrence relations is with respect to the number of main blocks.

\noindent\textbf{Theorem 4.} \textit{ For a fixed $j\in[n]$ the following formula is true
:}
\begin{equation}\label{r3}F(n,k,P,m)=\sum_{i=1}^{p_j}{p_j\choose i}F(n-1,k-i+1,P\setminus\{p_j\},m),\end{equation}

\noindent\textbf{Proof.} Omitting  $p_j$-block
 of $X$ we
obtain a  set $Y$ with $n-1$ main blocks, and additional $m$-block.
Each $n+k$-inset of $X$ may be obtained
as a union of some $n+k-i$- insets $(1\leq i\leq p)$ of $Y$ and
some of ${p_j\choose i}$-subsets of the omitting $m$-block, which proves (\ref{r3}).

In the case $p_1=p_2=\cdots=p_n=2,$ we shall denote $F(n,k,P,m)=c(n,k,m).$

In this case the  formula (\ref{r3}) gives
\begin{equation}\label{p2}c(n,k,p)=2c(n-1,k,m)+c(n-1,k-1,m).\end{equation}
We shall prove that this is the recurrence relation for coefficients of Chebyshev polynomials.

\noindent\textbf{Theorem 5.} \textit{ The numbers
$(-1)^kc(n,k,m),\;(m=0,1)$ are the coefficients of
$x^{n-k+m}$ in Chebyshev polynomial $P_{n+k+m}(x).$
If $m=0$ we obtain
 coefficients for the polynomials of the second, and for $m=1$ of the first kind.
\\In this way we may obtain almost all coefficients of Chebyshev polynomials of both kind.}

\noindent\textbf{Proof.}
We deal here with polynomials whose degree are $n+k+m,\;n\geq 1,\;k\geq 0.$ Thus, if $m=0$ we start with polynomials of the second  degree,
and with its $a(2,0)$ coefficient. If $m=1$ then we start with
polynomial of the third degree, and its $a(3,1)$ coefficient.

Multiplying (\ref{p2}) by $(-1)^k$ we obtain
\begin{equation}\label{cp}(-1)^kc(n,k,p)=2(-1)^kc(n-1,k,m)-(-1)^{k-1}c(n-1,k-1,m).\end{equation}
Denote $(-1)^kc(n,k,m)=a(r,s),$ where
\begin{equation}\label{sis}\begin{array}{c}n+k+m=r,\\n-k+m=s\end{array}.\end{equation}
This system has the solution for any $n,\;k$ if and only if $r$ and $s$ are of the same parity. This will be enough since only such Chebyshev coefficients
 are not zero.
The equation (\ref{cp}) becomes
$$a(n+k+m,n-k+m)=2a(n-1+k+m,n-1-k+m)-a(n-2+k+m,n-k+m),$$
that is $$a(r,s)=2a(r-1,s-1)-a(r-2,s),$$
which is well-known recurrence relation for Chebyshev coefficients.

For $m=0$ the initial condition is $a(2,0)=4$, which is the  coefficient
 of Chebyshev polynomials of the second kind, while in the case $m=1$
the initial condition is $a(3,1)=-3$  which is the coefficient of the polynomials of the first kind.

\noindent\textbf{Remark 1.} \textit{The system (\ref{sis}) has the solution by $n$ and $k$ if and only if $r$ and $s$ are of the same parity, which is the case
 with nonzero coefficients of Chebyshev polynomials. This means that in the preceding way we may obtain all Chebyshev coefficients except coefficients of polynomials
$U_0(x),T_0(x),T_1(x),T_2(x)$.}

According to (\ref{f1}) we also have

\noindent\textbf{Theorem 6.} \textit{If $a(r,s)$ is the coefficient of Chebyshev polynomial of degree $r$ of $x^s$ then
  $$(-1)^ka(r,k)=\sum_{i=0}^n(-1)^i{n\choose i}{2n+m-2i\choose
n+k}$$ where $n,\;k$ are the solution of the system (\ref{sis}). For $m=1$ we obtain coefficients for polynomial of the first kind, while $m=0$ produces the coefficients
of polynomials of the second kind.}

\noindent\textbf{Theorem 7.} \textit{ Let $I_0\not=\emptyset$ be arbitrary. Then
$$F(n,k,P,m)=F(n-r,k,Q,m)+\sum_{i=1}^{|I_0|}\sum_{j=0}^i{n\choose i}{i\choose j}F(n-r-j,k-i+j,Q,m).$$
where $r$ is the number of elements of $I_0$ such that $p_i=1,\;(i\in I_0),$
and $Q=\{q_1,q_2,\ldots,q_n\}$ where $q_i=p_i-1,\;(i\in I_0)$ and $q_i=p_i$ otherwise, and $q_i>0,$.}

\noindent\textbf{Proof.}
For nonempty $I_0\subseteq [n]$ take  $x_i$ from $p_i$-block, where  $i\in I_0$ and form the set $X(I_0)$ consisting of these elements. Denote further
 $Y=X\setminus X(I_0),$ and form the set $Q$ as in theorem.

Suppose first that the set  $\{i\in I_0: q_i=0\}$ is empty.
This means that $f(n,k,P,m)$ and $F(n,k,Q,m)$ have the same number of main blocks.

 Each $n+k$-inset of $Y$ is also an $n+k$-inset of $X.$ In such a way we obtain all $n+k$-insets of $X$ which do not meet $X(I_0).$ There are
$F(n,k,Q,m)$ such subsets.

If $U$ is some of remaining $n+k$-insets of $X,$ then $U$ intersects  $X(I_0).$ Take nonempty $I\subseteq I_0,$
and count those $n+k$-insets $U$ of $X$  such that $U\cap X(I_0)=X(I).$

First, such one inset may be obtained as a union of an $n+k-|I|$-inset of $Y$ with $X(I).$  There are $F(n,k-|I|,Q,m)$ such subsets.
Further, $U$ may be obtained as a union of some $n+k-|I|$ subset of $Y$ not intersecting  one, two, ...., $|I|$ blocks which indices belong to $I,$ with $X(I).$.
Thus, there are $$\sum_{j=0}^{|I|}F(n-j,k-|I|+j),$$ such subsets. Summing over all nonempty subsets of $I_0$ we finally obtain
$$F(n,k,P,m)=F(n,k,Q,m)+\sum_{i=1}^{|I_0|}\sum_{j=0}^i{n\choose i}{i\choose j}F(n-j,k-i+j,Q,m).$$

Suppose now that $\{i\in I_0: q_i=0\}$ is not empty. Denote $|\{i\in I_0: q_i=0\}|=r.$
This means that the set $Y$ has $n-r$ main blocks. The rest of the proof is the same as above.

\vspace{0.2cm}

\noindent\textbf{References}

\vspace{0.2cm}

\noindent[1] N. J. A. Sloane, (2007), The On-Line Encyclopedia of
Integer Sequences, published electronically at
www.research.att.com/~njas/sequences/.

\end{document}